# Some observations about Hankel determinants of the columns of Pascal's triangle and related topics

Johann Cigler

johann.cigler@univie.ac.at

**Abstract**

This note collects some facts and conjectures about the Hankel determinants and their generating functions of the columns of Hoggatt triangles which apparently are related to combinatorial objects such as Young tableaux and Narayana numbers.

## 0. Introduction

This note originated from the accidental observation that the Hankel determinants

$$d_k(m,r)=\det\left(\binom{k+i+j}{m}\right)_{i,j=0}^{r-1} \tag{1}$$

of the columns $\left(\binom{n}{m}\right)_{n\geq 0}$ of Pascal's triangle coincide (apart from the sign) with the columns $\left(\left\langle {k \atop m-r+1} \right\rangle_r\right)_{k\geq 0}$ of the $r-$Hoggatt triangle, an analog of Pascal's triangle, which has been introduced in [4] and further studied in [2] and [3]. The fact that the $r-$Hoggatt binomials $\left\langle {n \atop k} \right\rangle_r$ count the number of semistandard Young tableaux with shape $r^k$ and that their generating functions are related to Narayana polynomials and their $r-$dimensional analogs suggests that the Hankel determinants of the binomials $\binom{n}{m}$ also have some combinatorial interpretation related to these objects. It would be interesting to find such interpretations.

Computations suggest that also the Hankel determinants of the $r-$Hoggatt binomials

$$d_k(s,m,r)=\det\left(\left\langle {k+i+j \atop m} \right\rangle_s\right)_{i.j=0}^{r-1} \tag{2}$$

and their generating functions have similar properties. We formulate them as conjectures hoping that someone will prove them.

## 1. Hoggatt binomials

Let us first recall from [2] and [3] some properties of the $r-$Hoggatt binomials $\left\langle {n \atop k} \right\rangle_r$.

For a positive integer $r$ let

$$\langle n\rangle_r = \binom{n+r-1}{r} = \frac{(n+r-1)!}{r!(n-1)!} = \frac{(n)^{(r)}}{(1)^{(r)}} = \frac{(r+1)^{(n-1)}}{(1)^{(n-1)}} \tag{3}$$

and

$$\langle n\rangle_r ! = \langle 1\rangle_r \langle 2\rangle_r \cdots \langle n\rangle_r , \tag{4}$$

where for $r \in \mathbb{Z}$ the rising factorial $x^{(r)}$ is defined by $x^{(r)} = x(x+1)\cdots(x+r-1)$ if $r > 0,$ $x^{(0)} = 1$ and $x^{(r)} = \dfrac{1}{(x-1)(x-2)\cdots(x+r)}$ if $r < 0.$

The $r-$Hoggatt binomials $\left\langle {n \atop k} \right\rangle_r$ can be defined by

$$\left\langle {n \atop k} \right\rangle_r = \frac{\langle n\rangle_r !}{\langle k\rangle_r !\langle n-k\rangle_r !} = \prod_{j=1}^{k} \frac{\langle n-k+j\rangle_r}{\langle j\rangle_r} = \prod_{j=1}^{k} \frac{(n-k+j)^{(r)}}{(j)^{(r)}} = \prod_{j=1}^{k} \frac{(n-j+1)^{(r)}}{(k-j+1)^{(r)}} \tag{5}$$

for $0 \le k \le n$ and $\left\langle {n \atop k} \right\rangle_r = 0$ else.

Their generating functions can be expressed as

$$\sum_{k=0}^{n} \left\langle {n \atop k} \right\rangle_r t^k = {}_rF_{r-1}\left( {-n,-n-1,\cdots,-n-r+1 \atop 2,3,\cdots,r} ;(-1)^r t \right) \tag{6}$$

in terms of the hypergeometric function ${}_pF_q\left( {a_1,a_2,\cdots,a_p \atop b_1,b_2,\cdots,b_q} ;t \right) = \sum_{k\ge 0} \frac{(a_1)^{(k)}\cdots(a_p)^{(k)}}{(b_1)^{(k)}\cdots(b_q)^{(k)}} \frac{t^k}{k!}.$

This follows from

$$\left\langle {n \atop k} \right\rangle_r = \prod_{j=1}^{k} \frac{(n-k+j)^{(r)}}{(j)^{(r)}} = \prod_{j=1}^{k}\prod_{i=0}^{r-1} \frac{(n-k+j+i)}{(j+i)} = \frac{\prod_{j=0}^{r-1}(n-j+r-1)!\,j!}{\prod_{j=0}^{r-1}(n+j-k)!(k+r-j-1)!}$$

$$= \frac{\prod_{j=0}^{r-1}(n+j)!\,j!}{\prod_{j=0}^{r-1}(n+j-k)!(k+j)!} = \prod_{j=0}^{r-1} \frac{\binom{n+j}{k}}{\binom{k+j}{k}} = \prod_{j=0}^{r-1} \frac{(-n-j)^{(k)}}{(j+1)^{(k)}}(-1)^{rk}.$$

The last entry in (5) shows that the $r-$Hoggatt binomials can be interpreted as the number of semistandard Young tableaux with shape $r^k$ (a box with $r$ columns and $k$ rows) which generalizes the fact that $\binom{n}{k}$ is the number of $k-$element subsets of $\{1,2,\cdots,n\}$.

The formula $\binom{n}{k}=\frac{(n+1-k)^{(k)}}{(1)^k}$ extends to

$$\left\langle \begin{matrix} n \\ k \end{matrix} \right\rangle_r == \frac{\prod_{j=0}^{r-1}(n+1-k+j)^{(k+r-1-2j)}}{\prod_{j=0}^{r-1}(1+j)^{(k+r-1-2j)}}. \tag{7}$$

Noting that $\frac{(n+k)!}{n!}=(n+1)^{(k)}$ also holds for negative $k$ this follows from

$$\left\langle \begin{matrix} n \\ k \end{matrix} \right\rangle_r = \frac{\prod_{j=0}^{r-1}(n-j+r-1)!\,j!}{\prod_{j=0}^{r-1}(n+j-k)!(k+r-j-1)!} = \frac{\prod_{j=0}^{r-1}(n+1-k+j)^{(k+r-1-2j)}}{\prod_{j=0}^{r-1}(1+j)^{(k+r-1-2j)}}.$$

Setting

$$L_j = L_j(n,k,r) = \frac{(n+1-k+j)^{(k+r-1-2j)}}{(1+j)^{(k+r-1-2j)}} \tag{8}$$

(7) can be written as

$$\left\langle \begin{matrix} n \\ k \end{matrix} \right\rangle_r = L_0(n,k,r)\,L_1(n,k,r)\cdots L_{r-1}(n,k,r). \tag{9}$$

For example, for $k=1$ and $r=3$ we get

$$\langle n\rangle_3 = \left\langle \begin{matrix} n \\ 1 \end{matrix} \right\rangle_3 = L_0L_1L_2 = \frac{(n)^{(3)}}{(1)^{(3)}}\frac{(n+1)^{(1)}}{(2)^{(1)}}\frac{(n+2)^{(-1)}}{(3)^{(-1)}} = \frac{n(n+1)(n+2)}{3!}\frac{n+1}{2}\frac{2}{n+1} = \binom{n+2}{3}.$$

For each $n$ the entries $\left\langle \begin{matrix} n \\ k \end{matrix} \right\rangle_r$ are *palindromic with center of symmetry at* $\frac{n}{2}$,

$$\left\langle \begin{matrix} n \\ k \end{matrix} \right\rangle_r = \left\langle \begin{matrix} n \\ n-k \end{matrix} \right\rangle_r . \tag{10}$$

They are also *unimodal* with center of symmetry at $\frac{n}{2}$, which means that

$$\left\langle {n \atop 0} \right\rangle_r \le \left\langle {n \atop 1} \right\rangle_r \le \cdots \le \left\langle {n \atop \left\lfloor \frac{n}{2} \right\rfloor} \right\rangle_r = \left\langle {n \atop \left\lfloor \frac{n+1}{2} \right\rfloor} \right\rangle_r \ge \cdots \ge \left\langle {n \atop n-1} \right\rangle_r \ge \left\langle {n \atop n} \right\rangle_r .$$

In [3] it is shown that $\sum_{k=0}^{n} \left\langle {n \atop k} \right\rangle_r t^k$ is even gamma-positive. Recall that a polynomial $p(x)$ of degree $n$ is called gamma- positive if it can be written in the form $p(x) = \sum_{j=0}^{\left\lfloor \frac{n}{2} \right\rfloor} \gamma_j x^j (1+x)^{n-2j}$ with positive coefficients $\gamma_j$.

For $r = 2$ the numbers $\langle n \rangle_2 = \binom{n+1}{2} = T_n$ are the triangle numbers $1, 3, 6, 10, \cdots$, $\langle n \rangle_2 ! = \frac{n!(n+1)!}{2^n}$, and

$\left\langle {n \atop k} \right\rangle_2 = L_0(n,k,2) L_1(n,k,2) = \frac{(n+1-k)^{(k+1)}}{(1)^{k+1}} \frac{(n+2-k)^{(k-1)}}{(2)^{k-1}} = \frac{1}{k+1} \binom{n}{k} \binom{n+1}{k}$ are the Narayana numbers.

There is a nice generalization of the formula

$$\sum_{n\ge 0} \binom{n+k}{k} x^n = \frac{1}{(1-x)^{k+1}}. \tag{11}$$

In [6] Robert A. Sulanke introduced Narayana numbers $N(r,n,k)$ of dimension $r$. His results imply

$$(1-x)^{rs+1} \sum_k \left\langle {k+s \atop s} \right\rangle_r x^k = \sum_{j=0}^{(r-1)(s-1)} N(r,s,j) x^j. \tag{12}$$

For $r = 3$ the polynomials $\sum_{j=0}^{2(k-1)} N(3,k,j) x^j$ are $1$, $1 + 3x + x^2$, $1 + 10x + 20x^2 + 10x^3 + x^4$, $1 + 22x + 113x^2 + 119x^3 + 113x^4 + 22x^5 + x^6, \cdots$.

For $r = 2$ we get the Narayana numbers $N(2,n,k) = \frac{1}{k+1} \binom{n}{k} \binom{n-1}{k} = N_{n,k}$ in the usual notation. In our notation $N_{n,k} = \left\langle {n-1 \atop k} \right\rangle_2$.

**Remark**

There seems to be no known closed formula for the $r-$dimensional Narayana numbers, but their sum is the $r-$dimensional Catalan number

$$C_{r,n}=(rn)!\prod_{j=0}^{r-1}\frac{j!}{(n+j)!}. \tag{13}$$

## 2. Hankel determinants of the sequences $\left(\binom{n}{m}\right)_{n\geq 0}$.

For a polynomial $p(x)$ the Hankel determinants $\det\left(p(k+i+j)\right)_{i,j=0}^{n-1}$ vanish for $n>m+1$, because by elementary column operations the first column can be reduced to $\left(\Delta^{n-1}p(k+i)\right)_{i=0}^{n-1}=0$ since $\Delta^{m+1}p(x)=0.$

Let $p_m(x)=\binom{x}{m}$ for $m\in\mathbb{N}.$ We are interested in the Hankel determinants

$$d_k(m,r)=\det\left(\binom{k+i+j}{m}\right)_{i.j=0}^{r-1} \tag{14}$$

for fixed $r.$ Note that $d_k(m,r)=0$ for $k<m-r+1$ and $d_k\left(m-r+1,r\right)=(-1)^{\binom{r}{2}}.$

For $r=0$ it will be convenient to set $d_k(m,0)=1$ for $k\geq m$ and $d_k(m,0)=0$ else.

The case $r=1$ is trivial. But let us state some results which will later be generalized.

We have

$$d_k(m,1)=\binom{k}{m} \tag{15}$$

and

$$(1-x)^{m+1}\sum_{k\geq 0}d_k(m,1)x^k=(1-x)^{m+1}\sum_{k\geq 0}\binom{k}{m}x^k=x^m. \tag{16}$$

Computations suggest

**Theorem 1**

*The Hankel determinants* $d_k(m,r)$ *satisfy*

$$(-1)^{\binom{r}{2}}d_k(m+r-1,r)=\left\langle {k \atop m} \right\rangle_r \tag{17}$$

*for* $m\geq 0.$

Thus the matrix of the signed Hankel determinants is the $r-$Hoggatt matrix

$$\left((-1)^{\binom{r}{2}} d_k(m+r-1,r)\right)_{k,m\geq 0} = \left(\left\langle {k \atop j} \right\rangle_r\right)_{k,j\geq 0} \tag{18}$$

For example the sequences $\left(-d_k(2,3)\right)_{k\geq 0} = \left(1,1,1,\cdots\right)$, $\left(-d_k(3,3)\right)_{k\geq 0} = \left(0,1,4,10,20,35,\cdots\right)$, $\left(-d_k(4,3)\right)_{k\geq 0} = \left(0,0,1,10,50,175,490,\cdots\right)$, $\left(-d_k(5,3)\right)_{k\geq 0} = \left(0,0,0,1,20,175,980\cdots\right)$, $\cdots$ are the columns of the $3-$Hoggatt triangle

| $n \backslash k$ | 0 | 1 | 2 | 3 | 4 |
|---|---|---|---|---|---|
| 0 | 1 | | | | |
| 1 | 1 | 1 | | | |
| 2 | 1 | 4 | 1 | | |
| 3 | 1 | 10 | 10 | 1 | |
| 4 | 1 | 20 | 50 | 20 | 1 |
| 5 | 1 | 35 | 175 | 175 | 35 |
| 6 | 1 | 56 | 490 | 980 | 490 |
| 7 | 1 | 84 | 1176 | 4116 | 4116 |

**Proof of Theorem 1**

The condensation method for Hankel determinants (cf. [5]) gives for $m \geq r-1$

$$d_k\left(m,r\right)d_{k+2}\left(m,r-2\right)-d_{k+2}\left(m,r-1\right)d_k\left(m,r-1\right)+d_{k+1}\left(m,r-1\right)^2=0. \tag{19}$$

Let $D_k(m,r)=(-1)^{\binom{r}{2}}d_k(m,r).$ Then (19) is equivalent with

$$\frac{D_k\left(m,r\right)D_{k+2}\left(m,r-2\right)}{D_{k+1}\left(m,r-1\right)^2}+\frac{D_{k+2}\left(m,r-1\right)D_k\left(m,r-1\right)}{D_{k+1}\left(m,r-1\right)^2}=1. \tag{20}$$

Since (17) holds for $r=0$ and $r=1$ it suffices to prove (20) for $D_k(m,r)=\left\langle {k \atop m-r+1} \right\rangle_r$.

By (8) we have

$$D_k(m,r)=\left\langle {k \atop m-r+1} \right\rangle_r = \prod_{j=0}^{r-1}\frac{\left(k+r-m+j\right)^{(m-2j)}}{(1+j)^{(m-2j)}}, \tag{21}$$

$$D_{k+2}(m,r-2)=\left\langle {k+2 \atop m-r+3} \right\rangle_{r-2} = \prod_{j=0}^{r-3}\frac{\left(k+r-m+j\right)^{(m-2j)}}{(1+j)^{(m-2j)}}, \tag{22}$$

$$D_{k+1}(m,r-1)=\left\langle \begin{matrix} k+1 \\ m-r+2 \end{matrix} \right\rangle_{r-1}=\prod_{j=0}^{r-2}\frac{\left(k+r-m+j\right)^{(m-2j)}}{(1+j)^{(m-2j)}}, \tag{23}$$

$$D_{k}(m,r-1)=\left\langle \begin{matrix} k \\ m-r+2 \end{matrix} \right\rangle_{r-1}=\prod_{j=0}^{r-2}\frac{\left(k+r-m-1+j\right)^{(m-2j)}}{(1+j)^{(m-2j)}}, \tag{24}$$

$$D_{k+2}(m,r-1)=\left\langle \begin{matrix} k+2 \\ m-r+2 \end{matrix} \right\rangle_{r-1}=\prod_{j=0}^{r-2}\frac{\left(k+r-m+1+j\right)^{(m-2j)}}{(1+j)^{(m-2j)}}. \tag{25}$$

Setting

$$S_j=S_j\left(m,r\right)=\frac{\left(k+r-m+j\right)^{(m-2j)}}{(1+j)^{(m-2j)}} \tag{26}$$

(21), (22) and (23) give

$$\begin{aligned}\frac{D_k\left(m,r\right)D_{k+2}\left(m,r-2\right)}{D_{k+1}\left(m,r-1\right)^2}&=\frac{S_0\cdots S_{r-1}S_0\cdots S_{r-3}}{S_0^2\cdots S_{r-2}^2}=\frac{S_{r-1}}{S_{r-2}}\\&=\frac{\left(k+2r-m-1\right)^{(m-2r+2)}}{(r)^{(m-2r+2)}}\frac{(r-1)^{(m-2r+4)}}{\left(k+2r-m-2\right)^{(m-2r+4)}}=\frac{(r-1)(m-r+2)}{(k+1)\left(k+2r-m-2\right)}.\end{aligned} \tag{27}$$

On the other hand by (24), (25) and (23) we get

$$\begin{aligned}\frac{D_{k+2}\left(m,r-1\right)D_k\left(m,r-1\right)}{D_{k+1}\left(m,r-1\right)^2}&=\frac{\prod_{j=0}^{r-2}\frac{\left(k+r-m-1+j\right)^{(m-2j)}}{(1+j)^{(m-2j)}}\prod_{j=0}^{r-2}\frac{\left(k+r-m+1+j\right)^{(m-2j)}}{(1+j)^{(m-2j)}}}{\prod_{j=0}^{r-2}\frac{\left(k+r-m+j\right)^{(m-2j)}}{(1+j)^{(m-2j)}}\prod_{j=0}^{r-2}\frac{\left(k+r-m+j\right)^{(m-2j)}}{(1+j)^{(m-2j)}}}\\&=\prod_{j=0}^{r-2}\frac{\left(k+r-m-1+j\right)^{(m-2j)}\left(k+r-m+1+j\right)^{(m-2j)}}{\left(k+r-m+j\right)^{(m-2j)}\left(k+r-m+j\right)^{(m-2j)}}=\frac{(k-m+r-1)}{\left(k-m+2r-2\right)}\frac{(k+r)}{(k+1)}.\end{aligned} \tag{28}$$

By (27) and (28) we get

$$\frac{(r-1)(m-r+2)}{(k+1)\left(k+2r-m-2\right)}+\frac{(k-m+r-1)}{\left(k-m+2r-2\right)}\frac{(k+r)}{(k+1)}=1$$

which proves (20).

Consider

$$(1-x)^{r(m-r+1)+1}\sum_{k\ge 0} d_k(m,r)x^k = (-1)^{\binom{r}{2}}(1-x)^{r(m-r+1)+1}\sum_{k\ge 0}\left\langle \begin{matrix} k \\ m-r+1 \end{matrix}\right\rangle_r x^k.$$

From (12) we get

$$(1-x)^{r(m-r+1)+1}\sum_{k\ge 0}\left\langle \begin{matrix} k \\ m-r+1 \end{matrix}\right\rangle_r x^k = x^{m-r+1}(1-x)^{r(m-r+1)+1}\sum_{k\ge 0}\left\langle \begin{matrix} k+m-r+1 \\ m-r+1 \end{matrix}\right\rangle_r x^k$$

$$= x^{m-r+1}\sum_{j=0}^{(r-1)(m-r)} N(r,m-r+1,j)x^j.$$

Therefore, we get

**Theorem 2**

*The generating function of the Hankel determinants* $d_k(m,r)$ *satisfies*

$$(1-x)^{r(m-r+1)+1}\sum_{k\ge 0} d_k(m,r)x^k = (-1)^{\binom{r}{2}} x^{m-r+1}A_{m,r}(x) \tag{29}$$

*with* $A_{m,r}(x) = \sum_{j=0}^{(r-1)(m-r)} N(r,m-r+1,j)x^j.$

For example, $A_{m,1}(x)=1,\ A_{m,2}(x)=C_{m-1}(x)=\sum_{k=0}^{m-2} N_{m-1,k}x^k.$

For $r=3$ we have $A_{2,3}(x)=A_{3,3}(x)=1,\ A_{4,3}(x)=1+3x+x^2,$
$A_{5,3}(x)=1+10x+20x^2+10x^3+x^4.$

**Remark**

In J. Agapito [1] a proof is given that the polynomials $A_{m,r}(x)$ are unimodal and palindromic and even gamma-positive.

It is also shown that if

$$F_n(a,b,x)=\left(x^a D^b\right)^n \frac{1}{1-x} \tag{30}$$

then

$$A_{m,r}(x)=\frac{(1-x)^{r(m-r+1)+1}F_{m-r+1}\left(r-1,r,x\right)}{x^{r-1}\prod_{j=0}^{r-1}\frac{\left(m+j+1-r\right)!}{j!}} \tag{31}$$

for $m\ge r$ and $A_{r-1,r}(x)=1.$

## 3. Conjectures for Hankel determinants of the sequences $\left(\left\langle {n \atop m} \right\rangle_s\right)_{n\ge 0}$.

Computer experiments suggest analogous results for the Hankel determinants

$$d_k(s,m,r) = \det\left(\left\langle {k+i+j \atop m} \right\rangle_s\right)_{i.j=0}^{r-1}. \tag{32}$$

Theorem 1 implies that $\det\left(p(k+i+j)\right)_{i,j=0}^{r-1}$ is a polynomial in $k$ of degree $r(m-r+1)$ if $\deg p(x) = m.$ Since $\left\langle {n \atop m} \right\rangle_s$ has degree $ms$ we know that $d_k(s,m,r)$ has degree $r(ms-r+1).$

### Conjecture 3

$$(1-x)^{rsm-r^2+r+1}\sum_{k\ge 0} d_k(s,m,r)x^k = (-1)^{\binom{r}{2}} x^{m-r+1} A_{s,m,r}(x) \tag{33}$$

*where* $A_{s,m,r}(x)$ *is a polynomial of degree* $(rs-1)m-s-r^2+r+1$ *with positive coefficients which is gamma- positive. Further we have*

$$A_{s,m,r}(1) = C_{m,s}^r C_{r,ms-r+1}. \tag{34}$$

### Remark

It would be interesting if there exists an alternative formula for $A_{s,m,r}(x)$ which generalizes formula (31) or if there exists a combinatorial interpretation of identity (34).

### Conjecture 4

*Let*

$$S_j = S_j(s,m,r) = \frac{\left(k+r-m+j\right)^{(m+s-1-2j)}}{\left(1+j\right)^{(m+s-1-2j)}} \tag{35}$$

*and define*

$$w_k(s,m,r) = \left(S_0S_1\cdots S_{r-s}\right)^s \left(S_{r-s+1}S_{r-s+2}\right)^{s-1}\left(S_{r-s+3}S_{r-s+4}\right)^{s-2}\cdots\left(S_{r+s-5}S_{r+s-4}\right)^2\left(S_{r+s-3}S_{r+s-2}\right) \tag{36}$$

*Then we get for* $s \le r$

$$(-1)^{\binom{r}{2}} d_k(s,m,r) = \frac{w_k(s,m,r)u_k(s,m,r)}{w_{m-r+1}(s,m,r)u_{m-r+1}(s,m,r)} \tag{37}$$

$$(-1)^{\binom{s}{2}} d_k(r,m,s) = \frac{w_k(s,m+r-s,r)U_k(r,m,s)}{w_{m-s+1}(s,m+s-r,r)U_{m-s+1}(r,m,s)} \tag{38}$$

*where* $u_k(s,m,r)$ *and* $U_k(s,m,r)$ *are polynomials in* $\mathbb{Z}[x]$ *with degree*

$$\deg\left(u_k(s,m,r)\right) = (s-1)r^2 - \left(s^2-1\right)r + 2\binom{s+1}{3} \tag{39}$$

*and*

$$\deg U_k(r,m,s) = 2\binom{s}{3} \tag{40}$$

for $m \geq s-1$.

**Remark**

For $s=1$ formula (38) reduces to (9) and formula (37) to (17).

For $s=2$ we get

$$d_k(r,m,2) = -\left(\prod_{j=0}^{r-2} \frac{\left(k+2-m+j\right)^{(m+r-1-2j)}}{(1+j)^{(m+r-1-2j)}}\right)^2 \frac{\left\langle {k \atop m-r} \right\rangle_2}{\left\langle {m-1 \atop r-1} \right\rangle_2} = \frac{\left\langle {k+2 \atop m+1} \right\rangle_{r-1}^2 \left\langle {k \atop m-r} \right\rangle_2}{\left\langle {m-1 \atop m-r} \right\rangle_2}$$

For $s=3$ the polynomials $U_k(r,m,3)$ are quadratic polynomials

$$U_k\left(r,m,3\right) = C_{m,r}\left(k^2 + (r+4-m)k + \frac{\frac{3-(m-2)^2}{2}r^2 - (m-2)(2m-1)r - \frac{(m-2)^2+3}{2}}{mr-1}\right)$$